\documentclass[a4paper,12pt]{amsart}
\usepackage[letterpaper, margin=1.2in]{geometry}
\usepackage{latexsym}
\usepackage{amssymb}
\usepackage{amsmath}
\usepackage{amsfonts}
\usepackage[table]{xcolor}
\usepackage{pgfplots}
\usepackage{color}

\usepackage{txfonts,pxfonts} 
\usepackage{tikz}
\usetikzlibrary{calc,arrows}
\usetikzlibrary{calc}
\usepackage{verbatim}
\newtheorem{thm}{Theorem}[section]

\newtheorem{lem}[thm]{Lemma}

\newtheorem{cor}[thm]{Corollary}

\theoremstyle{definition}

\newtheorem*{rem}{Remark}

\def\fph{\mathbb{F}_{\ph}}

\newcommand{\Z}{\mathbb Z}
\newcommand{\z}{\mathbb Z}
\newcommand{\Q}{\mathbb Q}

\newcommand{\F}{\mathbb F}
\newcommand{\fp}{\mathbb F_p}

\def\F{\mathbb{F}}
\DeclareMathAlphabet{\mathpzc}{OT1}{pzc}{m}{it}

\newcommand{\p}{\mathfrak{p}}
\def\ol{\overline}
\def\al{\alpha}
\def\la{\lambda}

\def\th{\theta}
\def\md#1{\ \mbox{\rm(mod }{#1})}
\def\nph#1{N_{\ph}(#1)}
\def\npp#1{N_{\ph}^+(#1)}
\def\ph{\phi}
\newcounter{cs}
\stepcounter{cs}
\newcommand{\casos}{\begin{itemize}}
\newcommand{\fcasos}{\end{itemize}\setcounter{cs}{1}}

\newfont{\tit}{cmr12 scaled \magstep3}
\begin{document}
\title[]{On monogenity of certain  pure number fields defined by $x^{2^u\cdot 3^v\cdot 5^t}-m$}
\textcolor[rgb]{1.00,0.00,0.00}{}
\author{Lhoussain El Fadil}
\address{Faculty of Sciences Dhar El Mahraz, P.O. Box  1796 Atlas-Fez , Sidi Mohamed Ben Abdellah University, Fez-- Morocco}
\email{ lhoussain.elfadil@usmba.ac.ma}
\date{\today}
\maketitle
\begin{abstract}  
 Let $K = \mathbb{Q} (\alpha) $ be a pure number field generated by   a root $\alpha$ of a monic irreducible polynomial  $ F(x) = x^{2^u\cdot 3^v\cdot 5^t}-m$, with $ m \neq \pm 1 $   a square free rational integer,  $u$, $v$ and $t$  three  positive  integers.   In this paper, we study the monogenity of $K$. We prove that if $m\not\equiv 1\md4$, $m\not\equiv \pm 1\md9$, and $m\not\in\{\pm 1, \pm 7\}\md{25}$, then $K$ is monogenic. But if {$m\equiv 1\md{4}$} or $m\equiv 1\md9$ or $m\equiv -1\md9$ and $u=2k$ for some odd integer $k$ or $u\ge 2$ and $m\equiv 1\md{25}$ or $m\equiv -1\md{25}$ and $u=2k$ for some odd integer $k$ or $u=v=1$ and $m\equiv \pm 82\md{5^4}$, then $K$ is not monogenic.  
\end{abstract}
\section{Introduction} 
Let $K=\Q(\al)$ be a number field generated by    a  root { $\al$ } of a monic irreducible polynomial $F(x)\in
  \Z[x]$ and $\Z_K$  its ring of integers.
  It is well know that  the ring $\Z_K$ is a free $\Z$-module of rank $n=[K:\Q]$. Thus  by the fundamental theorem of finitely generated Abelian groups, the Abelian group $\Z_K/\Z[\al]$ is finite. Its cardinal  order is called the index of $\Z[\al]$, and denoted by $(\Z_K: \Z[\al])$.
  The ring $\Z_K$ is said to  be monogenic if it has a power integral basis  as a $\Z$-module. That is  $(1,\th,\cdots,\th^{n-1})$ is a $\Z$-basis of $\Z_K$ for some $\th\in \Z_K$. $K$ is said to be  not monogenic otherwise. 
  Monogenity of number fields is a classical problem of algebraic number theory, going back to Dedekind, Hasse, and Hensel (see for instance \cite{G19}). 
The problem of testing the monogenity of number fields and the construction of power integral bases has been intensively studied these last four decades, mainly by  Ga\'al, Nakahara,  Pohst, and their collaborators  (see for instance \cite{ AN,G, G19, GO, P}). 
 { In \cite{F4}, Funakura, calculated  integral bases and studied  monogenity of pure quartic fields.  In \cite{GR4},  Ga\'al and  Remete, calculated the elements of index $1$  in  pure quartic fields generated by $m^{\frac{1}{4}}$ for {$1< m <10^7$}  and $m\equiv 2,3 \md4$.  In \cite{AN6}, Ahmad, Nakahara, and Husnine proved  that  if $m\equiv 2,3 \md4$ and  $m\not\equiv \pm1\md9$, then the sextic number field generated by $m^{\frac{1}{6}}$ is monogenic.
They also showed in \cite{AN},    that if $m\equiv 1 \md4$ and $m\not\equiv \pm1\md9$, then the sextic number field generated by $m^{\frac{1}{6}}$ is not monogenic.  In \cite{E6}, based on prime ideal factorization, El Fadil showed that  if $m\equiv 1 \md4$ or $m\equiv 1\md9$, then the sextic number field generated by $m^{\frac{1}{6}}$ is not monogenic.
{Hameed and Nakahara proved that if $m\equiv 1\md4$, then the
octic number field generated by $m^{1/8}$ is not monogenic, but if $m\equiv 2,3 \md4$, 
then it is monogenic (\cite{HN8}).}   In \cite{GR17}, by applying the  explicit form of the index equation,{ Ga\'al and  Remete} obtained deep new   results on  monogenity of  number  fields generated  by $m^{\frac{1}{n}}$, with $3\le n\le 9$ and $m\neq \pm 1$ a square free integer. } In  \cite{Epr, E2r5s, E2u3v, E6, E6s,  E12, E24, E18,   E23k}, based on Newton polygon's techniques, El Fadil et al. studied the monogenity of some pure  number fields. Also  El Fadil, Chouli, and Kchit studied the monogenity of pure number fields defined by $x^{60}-m\in \Z[x]$ with $m$ a square free integer.
 The goal of this paper is to study the monogenity of pure number fields  defined by  $x^{2^u\cdot 3^v\cdot 5^t}-m$, where $m\neq \pm1$ is a square free integer, $u$, $v$, and $t$ are three natural integers.  The cases { $uvt=0$ have been studied in \cite{Epr, CNS, E2r5s, E2u3v}}. Also the case $u=2$ and $t=v=1$ has been  studied by  in \cite{E60}. Our proofs are based on Newton's polygon techniques  and on index divisors as introduced by Hensel as follows: The index of a field $K$ is defined by $i(K)=gcd\{(\mathbb{Z}_K:\mathbb{Z}[\theta])\mid K=\mathbb{Q}(\theta) \mbox{ and } \theta\in \Z_K \}$. A rational prime $p$ dividing $i(K)$ is called a prime common index divisor of $K$. If $\mathbb{Z}_K$ has a power integral basis, then $i(K)=1$. Therefore a field having a prime common index divisor is not monogenic. 
\section{Main results}
Let $K$ be a number field generated by  a complex root {$\al$} of   a monic  irreducible polynomial $F(x)=x^{2^u\cdot 3^v\cdot 5^t}-m$, with $m\neq \pm 1$  a square free rational integer, $u$, $v$, and $t$  three positive integers.
\begin{thm}\label{pib}
The ring $\Z[\al]$ is the ring of integers of $K$
 if and only if  $m\not\equiv 1\md4$, $m\not\equiv \pm 1\md9$, and $m\not\in\{\pm 1, \pm 7\}\md{25}$.
\end{thm}
\begin{rem}
If $m\equiv 1\md4$, $m\not\equiv \pm 1\md9$, and $m\not\in\{\pm 1, \pm 7\}\md{25}$, then $\Z[\al]$ is not integrally closed. But in this case, Theorem \ref{pib} cannot decide on monogenity of $K$. The following theorems give an answer.
\end{rem}
\begin{thm}\label{npib}
\begin{enumerate}
\item
$m\equiv 1\md{4}$, then $2$ divides $i(K)$.
\item
 $m\equiv 1\md9$ or $m\equiv -1\md{9}$ and $u=2k$ for some odd integer $k$ or $u=1$, $m\equiv -1\md{27}$, and  $v\ge 3$,  then $3$ divides $i(K)$.
  \item
 $m\equiv \pm 1\md{25}$ and $u= 2k$ for some odd integer $k$ or
$u=1$, $m\equiv \pm 1\md{125}$, and $t\ge 2$ or $u=v=1$,
 $m\equiv \pm 82\md{5^4}$, and $t\ge 3$,  then $5$ divides $i(K)$. 
\end{enumerate}
In particular, if one of these following conditions holds, then $K$ is not monogenic.  
\end{thm}
\begin{cor}\label{nspe}
Let  $a\neq \pm 1$ be a square free rational integer, $u$, $v$ and $t$  three positive integers, and $s< {2^u\cdot 3^v\cdot 5^t}$ a positive integer, which is coprime to $30$. Then $F(x)=x^{2^u\cdot 3^v\cdot 5^t}-a^s$ is irreducible over $\Q$. Let $K$ be the number field generated by  a complex root { $\al$} of   the monic  irreducible polynomial $F(x)$.
\begin{enumerate}
\item
If  $a\not\equiv 1\md4$,  $m\not\equiv \pm 1\md9$,  and $m\not\in\{\pm 1, \pm 7\}\md{25}$ then $K$ is monogenic.
\item
If one of the following conditions holds
\begin{enumerate}
\item
$a\equiv 1\md{4}$,
\item
 $a\equiv 1\md9$,
 \item
 $a\equiv -1\md{9}$ and $u= 2k$ for some positive odd integer $k$,
 \item
 $a\equiv -1\md{27}$ and  $v\ge 3$,
  \item
 $a\equiv  1\md{25}$ and $u\ge 2$,
 \item
 $a\equiv  -1\md{25}$ and  $u= 2k$ for some positive odd integer $k$,
 \item
 $a\equiv \pm 82\md{5^4}$, $u=v=1$, and $t\ge 3$,
\end{enumerate}
 then $K$ is not monogenic. 
 \end{enumerate}
\end{cor}
\section{Preliminaries}
{Let $K=\Q(\al)$ be a  number field generated by   a  root $\al$  of a monic irreducible polynomial $F(x)\in \Z[x]$, $\Z_K$ its ring of integers, and  $ind(\al)=(\Z_K: \Z[\al])$ the index of $\Z[\al]$ in $\Z_K$.
  For a rational prime integer $p$, if $p$ does not divide $(\mathbb{Z}_K : \mathbb{Z}[\alpha])$, then a well known theorem of Dedekind says that the  factorization of  $ p \mathbb{Z}_K$ can be derived directly  from the factorization of $\overline{F(x)}$ in $\F_p[x]$. Namely,  $p\Z_K=\prod_{i=0}^r \p_i^{l_i}$, where every $\p_i=p\Z_K+\phi_i(\al)\Z_K$ and $\ol{F(x)}=\prod_{i=1}^r \ol{\ph_i(x)}^{l_i}$ modulo $p$ is the factorization of $\ol{F(x)}$ into powers of monic irreducible coprime polynomials of $\F_p[x]$. So, $f(\p_i)={\mbox{deg}}(\phi_i)$ is the residue degree of $\p_i$ (see \cite[ Chapter I, Proposition 8.3]{Neu}).
  In order to apply this  theorem in an effective way, one needs a
criterion to test  whether  $p$ divides   the index $(\Z_K:\Z[\al])$. In $1878$,  Dedekind proved the well known Dedekind's criterion which allows to test if a prime integer $p$ divides $(\Z_K:\Z[\al])$ \cite[Theorem 6.1.4]{Co} and \cite{R}.
When Dedekind's criterion fails, that is,  $p$ divides the index $(\z_K:\z[\th])$  for every primitive element $\th\in \Z_K$ of $K$, {  then} it is not possible to obtain the prime ideal factorization of $p\z_K$ by applying { Dedekind's theorem}.
In 1928, Ore developed   an alternative approach
for obtaining the index $(\z_K:\z[\alpha])$, the
absolute discriminant, and the prime ideal factorization of the rational primes in
a number field $K$ by using Newton polygons (see \cite{MN, O}).
For the convenience of the reader, as it is necessary for the proof of our main results,  we refer to our paper \cite{GF8}.\\
  We start by recalling  some fundamental facts about Newton polygons  applied in algebraic number theory. For more details, we refer to \cite{El, EMN, GMN}. 
 For a  prime integer  $p$ and for a  monic polynomial 
$\phi\in \z[x]$ {whose reduction} is irreducible  in
$\fp[x]$, 
let $\fph$ be 
the field $\frac{\fp[x]}{(\overline{\phi})}$. For any
monic polynomial  $F(x)\in \z[x]$, upon  the Euclidean division
 by successive powers of $\ph$, we  expand $F(x)$ as
$F(x)=\sum_{i=0}^la_i(x)\phi(x)^{i}$, called    the $\phi$-expansion of $F(x)$
 (for every $i$, $deg(a_i(x))<
deg(\phi)$). 
The $\ph$-Newton polygon of $F(x)$ with respect to $p$, is the lower boundary of the convex envelope of the set of points $\{(i,\nu_p(a_i(x))),\, a_i(x)\neq 0\}$ in the Euclidean plane, which is denoted by $\nph{F}$. Let $S_1,\, S_2,\dots, S_t$ be the sides of $\nph{F}$. For every side $S$ of $\nph{F}$, the length of $S$, denoted by $l(S)$, is the length of its projection to the $x$-axis, its height, denoted by $H(S)$, is the length of its projection to the $y$-axis. Let $\la=H(S)/l(S)$, then $-\la$ is the slope of $S$. If  $\la\neq 0$, then  $\la=h/e$ with $e$ and $h$  two positive coprime integer.  Notice that $e=l(S)/d$, called the ramification index of $S$ and $h=H(S)/d$, where   $d=$gcd$(l(S),H(S))$ is called the degree of $S$. Thus $\nph{F}$ is the join of its different sides ordered by increasing slopes, which we can express by $\nph{F}=S_1+ S_2+\dots+ S_t$.
 The principal $\ph$-Newton polygon of $F(x)$ ,
 denoted by $\npp{F}$, is the part of the  polygon $\nph{F}$, which is  determined by  all sides of negative  slopes of $\nph{F}$.
  For every side $S$ of {$\npp{F}$}, with initial point $(s, u_s)$ and length $l$, and for every 
$i=0, \dots,l$, we attach   the following
 residue coefficient $c_i\in\fph$ as follows:
 $$c_{i}=
\left
\{\begin{array}{ll} 0,& \mbox{ if } (s+i,{\it u_{s+i}}) \mbox{ lies strictly
above } S
 ,\\
\left(\dfrac{a_{s+i}(x)}{p^{{\it u_{s+i}}}}\right)
\,\,
\md{(p,\phi(x))},&\mbox{ if }(s+i,{\it u_{s+i}}) \mbox{ lies on }S.
\end{array}
\right.$$
where $(p,\phi(x))$ is the maximal ideal of $\z[x]$ generated by $p$ and $\ph$.\\
Let  $\la=-h/e$ be the slope of $S$, where  $h$ and $e$ are two positive coprime integers. Then $d=l/e$ is the degree of $S$.  Notice that, the points  with integer coordinates lying on $S$ are exactly $(s,u_s),(s+e,u_{s}-h),\cdots, (s+de,u_{s}-dh)$. Thus, if $i$ is not a multiple of $e$, then 
$(s+i, u_{s+i})$ does not lie on $S$, and so $c_i=0$. Let
$F_S(y)=t_dy^d+t_{d-1}y^{d-1}+\cdots+t_{1}y+t_{0}\in\fph[y]$, called  
the residual polynomial of $F(x)$ associated to the side $S$, where for every $i=0,\dots,d$,
 $t_i=c_{ie}$.
The  theorem of Ore plays a key role for proving our main theorems:\\
  Let $\ph\in\Z[x]$ be a monic polynomial, with $\overline{\ph(x)}$  irreducible in $\F_p[x]$. As defined in \cite[Def. 1.3]{EMN},   the $\ph$-index of $F(x)$, denoted  $ind_{\ph}(F)$, is  deg$(\ph)$ multiplied by  the number of points with natural integer coordinates that lie below or on the polygon $\npp{F}$, strictly above the horizontal axis, {and strictly beyond the vertical axis} (see Figure 1).
 \begin{figure}[htbp] 
\begin{tikzpicture}[x=1cm,y=0.5cm]
\draw[latex-latex] (0,6) -- (0,0) -- (10,0) ;
\draw[thick] (0,0) -- (-0.5,0);
\draw[thick] (0,0) -- (0,-0.5);
\node at (0,0) [below left,blue]{\footnotesize $0$};
\draw[thick] plot coordinates{(0,5) (1,3) (5,1) (9,0)};
\draw[thick, only marks, mark=x] plot coordinates{(1,1) (1,2) (1,3) (2,1)(2,2)     (3,1)  (3,2)  (4,1)(5,1)  };
\node at (0.5,4.2) [above  ,blue]{\footnotesize $S_{1}$};
\node at (3,2.2) [above   ,blue]{\footnotesize $S_{2}$};
\node at (7,0.5) [above   ,blue]{\footnotesize $S_{3}$};
\end{tikzpicture}
\caption{    \large  $\npp{F}$.}
\end{figure}
 Assume that $\overline{F(x)}=\prod_{i=1}^r\overline{\ph_i}^{l_i}$ is the factorization of $\overline{F(x)}$ in $\F_p[x]$,  {where  $\ph_1,\dots, \ph_r$ are monic polynomials lying in $\Z[x]$ and $\ol{\ph_1},\dots, \ol{\ph_r}$ are pairwise coprime irreducible polynomials over $\F_p$}.
For every $i=1,\dots,r$, let  $N_{\ph_i}^+(F)=S_{i1}+\dots+S_{ir_i}$ be the principal part of the $\ph_i$-Newton polygon of $F$ with respect to $p$. For every {$j=1,\dots, r_i$},  let $F_{S_{ij}}(y)=\prod_{s=1}^{s_{ij}}\psi_{ijs}^{a_{ijs}}(y)$ be the factorization of $F_{S_{ij}}(y)$ { into powers of monic irreducible polynomials } of $\F_{\ph_i}[y]$. 
  Then we have the following  theorem of Ore (see \cite[Theorem 1.7 and Theorem 1.9]{EMN}, \cite[Theorem 3.9]{El}, and \cite{MN}):
 \begin{thm}\label{ore} (Theorem of Ore)
 \begin{enumerate}
 \item
  {$$\nu_p((\z_K:\z[\al]))\ge \sum_{i=1}^r ind_{\ph_i}(F).$$} 
  The equality holds if $a_{ijs}=1$ for every $i,j,s$.
\item
If  $a_{ijs}=1$ for every $i,j,s$, then
$$p\Z_K=\prod_{i=1}^r\prod_{j=1}^{r_i}
\prod_{s=1}^{s_{ij}}\p^{e_{ij}}_{ijs}$$ is the factorization of $p\Z_K$ into powers of prime ideals of $\Z_K$ lying above $p$, where $e_{ij}$ is the ramification index
 of the side $S_{ij}$ and $f_{ijs}=\mbox{deg}(\ph_i)\times \mbox{deg}(\psi_{ijs})$ is the residue degree of $\p_{ijs}$ over $p$ for every $i=1,\dots,r$, $j=1,\dots,r_i$, and $s=1,\dots, s_{ij}$.
 \end{enumerate}
\end{thm}
\begin{cor}\label{indore}
{ Under the assumptions above Theorem \ref{ore}, if for every $i=1,\dots,r$, $l_i=1$ or $N_{\ph_i}(F)=S_i$ has a single side of height $1$, then $\nu_p((\z_K:\z[\al]))=0$}.
\end{cor}
The following lemma allows to evaluate  the $p$-adic valuation of the binomial coefficient $\binom{p^r}{j} $. For its proof, we refer  to \cite{E2r5s}.
\begin{lem} \label{binomial}
	Let $p$ be a rational prime integer and $r$ be a positive integer. Then $ \nu_p\left(\binom{p^r}{j}\right)  =  r - \nu_p(j)$
	for any integer $j= 1,\dots,p^r-1 $. 
\end{lem}
The following lemma  allows to determine the $\ph$-Newton polygon of $F(x)$.
\begin{lem}\label{NP}
Let $F(x)=x^n-m\in \Z[x]$ be an irreducible polynomial and $p$ a prime integer which divides $n$ and does not divide $m$. Let $n=p^rt$ in $\Z$ with $p$ does not divide $t$. Then $\ol{F(x)}=\ol{(x^t-m)}^{p^r}$. Let $v=\nu_p(m^p-m)$ and $\ph\in \Z[x]$ be a monic polynomial, whose reduction modulo $p$ divides $\ol{F(x)}$. 
 Let us denote $(x^t-m)=\ph(x)Q(x)+R(x)$. Then $\nu_p(R)\ge1$. 
 \begin{enumerate}
 \item
 If $\nu_p(m^{p-1}-1)\le r$, then $\npp{F}$ is the lower boundary of the  convex envelope of the set of the points $\{(0,v)\}\cup \{(p^j,r-j), \, j=0,\dots,r\}$.
 \item
 If $\nu_p(m^{p-1}-1)\ge r+1$, then $\npp{F}$ is the lower boundary of the  convex envelope of the set of the {points} $\{(0,V)\}\cup \{(p^j,r-j), \, j=0,\dots,r\}$ for some integer $V\ge r+1$.
 \end{enumerate}
\end{lem}}
\section{Proofs of main results}
\begin{proof} of Theorem \ref{pib}.\\
The proof of Theorem \ref{pib} can be done by using Dedekind's criterion as it was shown in the proof of \cite[Theorem $6.1$]{HS}. But as the other results are based on Newton polygon's techniques, let us use theorem of index with "if and only if" as it is given in  \cite[Theorem 4.18]{GMN}, which  says that: $\nu_p(\Z_K:\Z[\al])=0$ if and only if $ind_1(F)=0$, where $ind_1(F)$ is the index given in Theorem \ref{ore}.  Since $\triangle(F)  =\mp (2^{u}\cdot 3^{v}\cdot 5^{t})^{2^{u}\cdot 3^{v}\cdot 5^{t}} m ^{2^{u}\cdot 3^{v}\cdot 5^{t}-1}$,  by the  formula  $\nu_p(\triangle(F))=2\nu_p(ind(F))+\nu_p(d_K)$, where $d_K$ is the absolute discriminant of $K$ and $ ind(F)=(\Z_K: \Z[\al])$, we conclude that $\mathbb{Z}[\alpha]$ is integrally closed if and only if   $p$ does not divide $(\Z_K: \Z[\al])$ for every rational prime integer $p$ dividing $30m$.  Let $p$ be a rational prime dividing $m$, then $ F(x) \equiv \ph^{2^{u}\cdot 3^{v}\cdot 5^{t}}  (\ mod \ p )$, where $\phi = x$. As $m$ is a square free integer, the $\phi$-principal Newton polygon with respect to $\nu_p$,  $ \npp {F} = S $ has a single side of height $\nu_p(m)$. As $l(S)=2^{u}\cdot 3^{v}\cdot 5^{t}\ge 2$, $ind_\ph(F)=0$ if and only if $S$ has  height $1$, which  means $\nu_p(m)=1$. 
  It follows that the unique prime candidates to divide the index $(\mathbb{Z}_K:\mathbb{Z}[\alpha])$ are $2$, $3$, and $5$.\\
 For $p=2$  and $2$ does not divide $m$,  
 let $\ph\in \Z[x]$ be a monic polynomial, whose reduction is an irreducible factor of  $(x^{3^{v}\cdot 5^{t}}-1)$ in $\F_2[x]$.
 Again  as $l(\npp{F})=2^{u}\ge 2$, $ind_\ph(F)=0$ if and only if $\npp{F}$ has a single side of  height  $1$, which  means by Lemma \ref{NP} that $ \nu_2(1- m)=1$; $m\not\equiv 1\md4$.\\  
 Similarly, for  $p=3$  and $3$ does not divide $m$,  let   $\ph\in \Z[x]$ be a monic polynomial, whose reduction is an irreducible factor of  $(x^{2^{u}\cdot 5^{t}}-m)$ in $\F_3[x]$.
 Again  as $l(\npp{F})=3^{v}\ge 2$, $ind_\ph{F}=0$ if and only if $\npp{F}$ has a single side of height $1$, which  means by Lemma \ref{NP} that $ \nu_3(m^2- 1)=1=1$; $m\not\equiv \mp1\md9$.  \\
 Again, for  $p=5$  and $5$ does not divide $m$,  let   $\ph\in \Z[x]$ be a monic polynomial, whose reduction is an irreducible factor of  $(x^{2^{u}\cdot 3^{v}}-m)$ in $\F_5[x]$.
 As $l(\npp{F})=5^{t}\ge 2$, $ind_\ph{F}=0$ if and only if $\npp{F}$ has a single side of height $1$, which  means by Lemma \ref{NP} that $ \nu_5(m^4- 1)=1=1$; $m\not\in \{\mp1, \pm 7\}\md{25}$.  
  \end{proof} 
 \smallskip
	
	{ The existence of prime common index divisors was first established in $1871$ by Dedekind who exhibited examples in fields of third and fourth degrees, for example, he considered the cubic field $K$ defined by $F(x)=x^3-x^2-2x-8$ and he showed that the prime $2$ splits completely. So, if we suppose that $K$ is monogenic, then we would be able to find a cubic polynomial generating $K$, that splits completely into distinct polynomials of degree $1$ in $\mathbb{F}_2[x]$. Since there are only $2$ distinct polynomials of degree $1$ in $\mathbb{F}_2[x]$, this is impossible. Based on these ideas and using Kronecker's theory of algebraic numbers, Hensel  gave a necessary and sufficient condition on the so-called "index divisors" for any prime integer  $p$ to be  a prime common index divisor \cite{He2}. $($For more details see \cite{HS}$)$}.
	{For the proof of Theorem $\ref{npib}$, we need the following lemma, which  characterizes the prime common index divisors of $K$.  We need to use only  one way, which is an immediate consequence of  Dedekind's theorem}.
	 	\begin{lem}\label{index}
		Let $p$ be a rational prime integer and $K$ be a number field. For every positive integer $f$, let $\mathcal{P}_f$ be the number of distinct prime ideals of $\mathbb{Z}_K$ lying above $p$ with residue degree $f$ and $\mathcal{N}_f$ the number of monic irreducible polynomials of $\mathbb{F}_p[x]$ of degree $f$.  Then $p$ is a prime common index divisor of $K$ if and only if $\mathcal{P}_f>\mathcal{N}_f$ for some positive integer $f$.
	\end{lem}
	\begin{rem}
		{As it was shown in the proof of Theorem \ref{pib}, the unique prime candidates to be a prime common index divisors of $K$ are $2$, $3$, and $3$, because if $p\not\in\{2, 3, 5\}$, then $p$ does not divide the index $(\mathbb{Z}_K:\mathbb{Z}[\alpha])$, and so the factorization of $p\mathbb{Z}_K$ is analogous to the factorization of $x^{2^u\cdot 3^v\cdot 5^t}-m$ in $\mathbb{F}_p[x]$.}
	\end{rem}
	\begin{rem}\label{remore}
		{In order to prove  Theorem \ref{npib}, we don't need to determine the factorization of $p\Z_K$ explicitly. But according to Lemma \ref{index}, we need only to show that $\mathcal{P}_f>\mathcal{N}_f$ for an adequate positive integer  $f$. So in practice the second point of Theorem \ref{ore}, could be replaced by the following:
			If  $l_i=1$ or $d_{ij}=1$ or $a_{ijk}=1$ for some $(i,j,k)$ according to notation of    Theroem \ref{ore}, then $\psi_{ijk}$ provides  a prime ideal $\p_{ijk}$ of $\Z_K$ lying above $p$ with residue degree  {$f_{ijk}=m_i\times t_{ijk}$}, where  $t_{ijk}=$deg$(\psi_{ijk})$ and $p\Z_K=\p_{ijk}^{e_ij}I$, where the factorization of the ideal $I$ can be derived from the other factors of each residual polynomials of $F(x)$.}
	\end{rem}
	 \begin{proof} of Theorem \ref{npib}.
 \begin{enumerate}
 \item
  Assume that   $m\equiv 1\md4$. Let $3^v\cdot 5^t=3s$ for some odd integer $s\in \Z$. Then $\overline{F(x)}=\overline{((x^{3}-1)(U(x))^{2^u}}=\overline{(x^2+x+1)T(x)}^{2^u}$  in $\F_2[x]$ for some monic polynomials $U$ and $T$ in $\Z[x]$  such that $\overline{(x^{2}+x+1)}$ and $\overline{T(x)}$ are coprime over $\F_2[x]$ because $\overline{(x^{3s}-1)}$ is separable over $\F_2[x]$. Let  $\phi={x^2+x+1}$ and $x^{3s}-1=  \ph(x) T(x)+2a$ for some integer $a\in \Z$. 
  \begin{enumerate}
 \item
{ $\nu_2(1-m)=2$, then by Lemmas \ref{NP} and \ref{binomial}, $\npp{F}=S$ has a single side of degree $d=2$. 
By using 
$F(x) = ((x^{3s}-1)+1)^{2^u}-m=(\ph(x) T(x)+2a)^{2^u}+\sum_{j=1}^{2^u-1}\binom{2^u}{j}(T\ph+2a)^{3^v - j}\cdots+1-m$, we have 
 $F_{S}(y)=t^2y^2+ty+1$, where $t\equiv T(x)\md{2,\ph}$ is a nonzero element of $\fph$ (because $\ol{\ph}$ does not divide $\ol{T(x)}$ in $\F_2[x]$). Hence $F_{S}(y)=(ty-x)(ty-x^2)$ in $\fph[y]$. Thus by Remark \ref{remore}, $\ph$ provides $2$ distinct prime ideals of $\Z_K$ lying above $2$ with residue degree $2$ each.}
 If $\nu_2(1-m)=3 $ , then by  Lemmas \ref{NP} and \ref{binomial}, $\npp{F}$ has  two sides $S_{1}$ and $S_{2}$ joining the point $(0,3)$, $(2^{u-1},1)$, and  $(2^u,0)$ (see $FIGURE\, 2$).    Thus $S_{1}$ is a side of  degree $2$  and $S_{2}$ is a side of  degree $1$. By using $F(x) = ((x^{3s}-1)+1)^{2^u}-m=(\ph(x) T(x)+2a)^{2^u}+\sum_{j=1}^{2^u-1}\binom{2^u}{j}(T\ph+2a)^{3^v - j}\cdots+1-m$, we conclude that 
  $F_{S_1}(y)=t^2y^2+ty+1=(ty-x)(ty-x^2)$, where $t\equiv T(x)\md{2,\ph}$ is a nonzero element of $\fph$  and  $F_{S_{2}}(y)$ is of degree $1$. By Remark \ref{remore},   $\ph$ provides three prime ideals of $\Z_K$ lying above $2$ with  residue degree deg$(\ph)=2$ each. As $x^2+x+1$ is the unique monic irreducible polynomial of degree $2$ in $\F_2[x]$,  by Lemma \ref{index}, $2$ divides $i(K)$ and $K$ is not monogenic. 
\begin{figure}[htbp] 
\begin{tikzpicture}[x=0.5cm,y=0.5cm]
\draw[latex-latex] (0,6.5) -- (0,0) -- (28,0) ;
\draw[thick] (0,0) -- (-0.5,0);
\draw[thick] (0,0) -- (0,-0.5); 
\draw[thick,red] (-2pt,1) -- (2pt,1);
\draw[thick,red] (-2pt,2) -- (2pt,2);
\draw[thick,red] (-2pt,3) -- (2pt,3);
\node at (0,0) [below left,blue]{\footnotesize  $0$};
\node at (6,0) [below ,blue]{\footnotesize  $2^{u-2}$};
\node at (13,0) [below ,blue]{\footnotesize  $2^{u-1}$};
\node at (26,0) [below ,blue]{\footnotesize  $2^{u}$};
\node at (0,1) [left ,blue]{\footnotesize  $1$};
\node at (0,2) [left ,blue]{\footnotesize  $2$};
\node at (0,3) [left ,blue]{\footnotesize  $3$};
\draw[thick,mark=*] plot coordinates{(0,3)(6.2,2)(13,1) (26,0)};
\node at (4,2) [above right  ,blue]{\footnotesize  $S_{1}$};
\node at (18,1) [above right  ,blue]{\footnotesize  $S_{2}$};
\end{tikzpicture}
\caption{    \large   $\npp{F}$ for  $v_2=3$.\hspace{5cm}}
\end{figure}
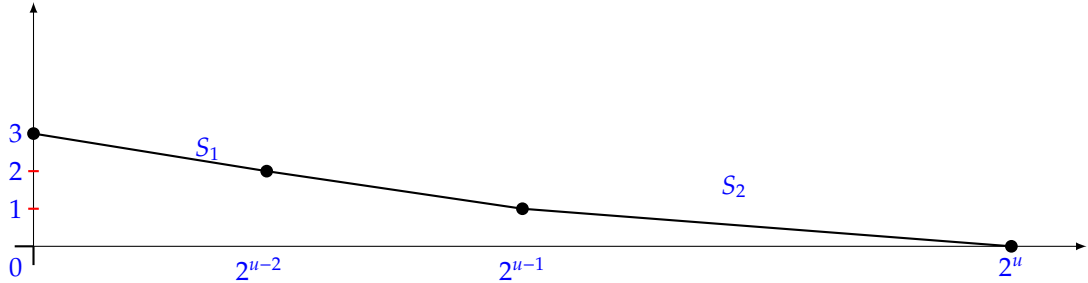
 \item
 If    $v_2\ge 4$, then by Lemma \ref{NP}, $\npp{F}$ has  $g\ge 2$ sides for which the last two sides $S_{g}$ and $S_{g-1}$ are of height $1$  each (see $FIGURE 3$).    Thus,  $F_{S_{g}}(y)$ and $F_{S_{g-1}}(y)$ are of degree $1$. By Remark \ref{remore},   $\ph$ provides at least two  prime ideals  of $\Z_K$ lying above $2$ with  residue degree deg$(\ph)=2$ each. As $x^2+x+1$ is the unique monic irreducible polynomial of degree $2$ in $\F_2[x]$,  by Lemma \ref{index}, $2$ divides $i(K)$ and $K$ is not monogenic. 
\begin{figure}[htbp] 
\begin{tikzpicture}[x=0.5cm,y=0.5cm]
\draw[latex-latex] (0,6.5) -- (0,0) -- (28,0) ;
\draw[thick] (0,0) -- (-0.5,0);
\draw[thick] (0,0) -- (0,-0.5); 
\draw[thick,red] (-2pt,1) -- (2pt,1);
\draw[thick,red] (-2pt,2) -- (2pt,2);
\draw[thick,red] (-2pt,3) -- (2pt,3);
\node at (0,0) [below left,blue]{\footnotesize  $0$};
\node at (6,0) [below ,blue]{\footnotesize  $2^{u-2}$};
\node at (13,0) [below ,blue]{\footnotesize  $2^{u-1}$};
\node at (26,0) [below ,blue]{\footnotesize  $2^{u}$};
\node at (0,1) [left ,blue]{\footnotesize  $1$};
\node at (0,2) [left ,blue]{\footnotesize  $2$};
\node at (0,3) [left ,blue]{\footnotesize  $3$};
\node at (0,4.5) [left ,blue]{\footnotesize  $\nu_2$};
\draw[thick,mark=*] plot coordinates{(6.2,2)(13,1) (26,0)};
\node at (8,2) [above right  ,blue]{\footnotesize  $S_{g-1}$};
\node at (18,1) [above right  ,blue]{\footnotesize  $S_{g}$};
\end{tikzpicture}
\caption{    \large   $\npp{F}$ for  $v_2\ge 4$.\hspace{5cm}}
\end{figure}
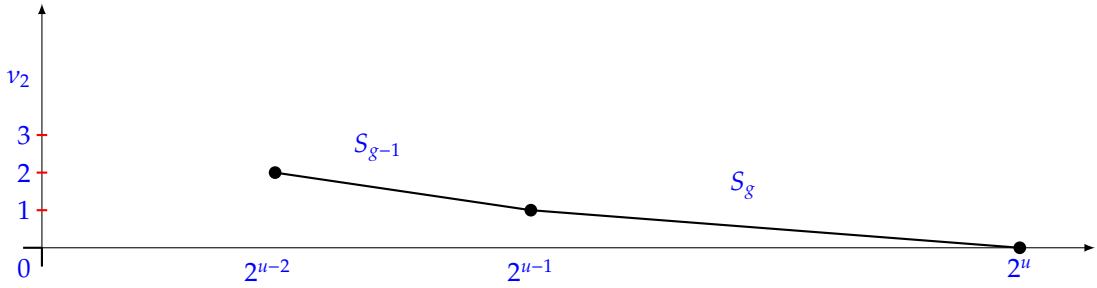
 \end{enumerate}
 \item
Assume  $m\equiv 1\md9$. Let $2^u\cdot 5^t=2s$ for some integer $s\in \Z$.
 Then $\overline{F(x)}=(x^{2s}-1)^{3^v}=((x-1)(x+1)U(x))^{3^v}$ in $\F_3[x]$ for some  monic polynomial  $U(x)\in \F_3[x]$. Let $\ph_1=x-1$, $\ph_2=x+1$, and $v_3=\nu_3(1-m)$.  If  $v_3\ge 2$, then by Lemmas \ref{NP} and \ref{binomial}, $N_{\ph_i}^+(F)$ has  $g\ge 2$ sides of which the last two sides $S_{ig}$ and $S_{ig-1}$ are of height $1$ each for every $i=1,2$ (see $FIGURE 4$ and $FIGURE 5$).
Thus  $F_{S_{ig}}(y)$ and $F_{S_{ig-1}}(y)$ are of degree $1$ for every $i=1,2$. By  Remark \ref{remore}, every  $\ph_i$ provides at least $2$  prime ideals $\p_{ij}$ of $\Z_K$ lying above $3$ with  residue degree $f_{ij}=1$ for every $i,j=1,2$. Therefore, there are at least $4$ prime ideals of $\Z_K$ lying above $3$ with  residue degree $1$ each. As  there is only $3$ monic irreducible polynomial of degree $1$ in $\F_3[x]$,  by Lemma \ref{index}, $3$ divides $i(K)$ and $K$ is not monogenic. 
\begin{figure}[htbp] 
\begin{tikzpicture}[x=0.5cm,y=0.5cm]
\draw[latex-latex] (0,6.5) -- (0,0) -- (28,0) ;
\draw[thick] (0,0) -- (-0.5,0);
\draw[thick] (0,0) -- (0,-0.5); 
\draw[thick,red] (1,-2pt) -- (1,2pt);
\draw[thick,red] (3,-2pt) -- (3,2pt);
\draw[thick,red] (9,-2pt) -- (9,2pt);
\draw[thick,red] (-2pt,1) -- (2pt,1);
\draw[thick,red] (-2pt,2) -- (2pt,2);
\draw[thick,red] (-2pt,3) -- (2pt,3);
\node at (0,0) [below left,blue]{\footnotesize  $0$};
\node at (3,0) [below ,blue]{\footnotesize $3^{v-2}$};
\node at (9,0) [below ,blue]{\footnotesize  $3^{v-1}$};
\node at (27,0) [below ,blue]{\footnotesize  $3^{v}$};
\node at (0,1) [left ,blue]{\footnotesize  $1$};
\node at (0,2) [left ,blue]{\footnotesize  $2$};
\node at (0,3) [left ,blue]{\footnotesize  $v_3$};
\draw[thick,mark=*] plot coordinates{(0,3) (3,2) (9,1) (27,0)};
\node at (6,1.3) [above right  ,blue]{\footnotesize  $S_{i2}$};
\node at (16,0.5) [above right  ,blue]{\footnotesize  $S_{i1}$};
\end{tikzpicture}
\caption{    \large   $N_{\ph_i}^+({F})$ for  $v\ge 2$ and $v_3=3$.\hspace{5cm}}
\end{figure}
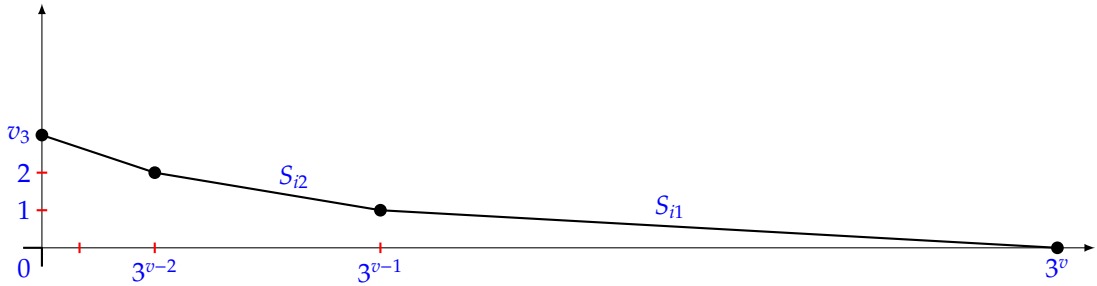
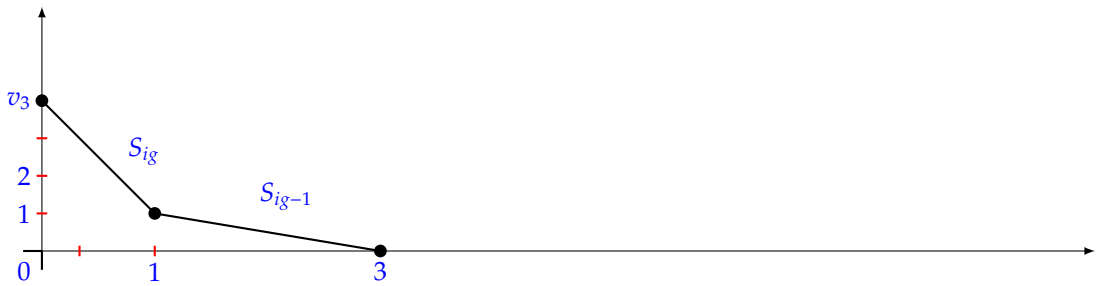
\begin{figure}[htbp] 
\begin{tikzpicture}[x=0.5cm,y=0.5cm]
\draw[latex-latex] (0,6.5) -- (0,0) -- (28,0) ;
\draw[thick] (0,0) -- (-0.5,0);
\draw[thick] (0,0) -- (0,-0.5); 
\draw[thick,red] (1,-2pt) -- (1,2pt);
\draw[thick,red] (3,-2pt) -- (3,2pt);
\draw[thick,red] (9,-2pt) -- (9,2pt);
\draw[thick,red] (-2pt,1) -- (2pt,1);
\draw[thick,red] (-2pt,2) -- (2pt,2);
\draw[thick,red] (-2pt,3) -- (2pt,3);
\node at (0,0) [below left,blue]{\footnotesize  $0$};
\node at (3,0) [below ,blue]{\footnotesize $1$};
\node at (9,0) [below ,blue]{\footnotesize  $3$};
\node at (0,1) [left ,blue]{\footnotesize  $1$};
\node at (0,2) [left ,blue]{\footnotesize  $2$};
\node at (0,4) [left ,blue]{\footnotesize  $v_3$};
\draw[thick,mark=*] plot coordinates{(0,4) (3,1) (9,0)};
\node at (2,2) [above right  ,blue]{\footnotesize  $S_{ig}$};
\node at (5.5,0.8) [above right  ,blue]{\footnotesize  $S_{ig-1}$};
\end{tikzpicture}
\caption{    \large   $N_{\ph_i}^+({F})$ for  $v=1$ and $v_3\ge 3$.\hspace{5cm}}
\end{figure}
\item
 Assume that  $m\equiv -1\md9$ and $u=2k$ for some odd integer $k$. Then $2^u\cdot 5^t=4s$ for some odd nonnegative integer $s$. Thus $\overline{F(x)}=\overline{(x^{4}+1)U(x)^{3^v}}=((x^{2}+x-1)(x^2-x-1)U(x))^{3^v}$ in $\F_3[x]$ for some monic polynomial $U\in\Z[x]$. Let $\ph_1=x^2+x-1$, $\ph_2=x^2-x-1$, and $v_3=\nu_3(1-m)$. Since $3$ does not divide $4s=2^u\cdot 5^t$, $\overline{(x^{4s}+1)}$ is separable over $\F_3$. So, $\overline{\ph_1\ph_2(x)}$ and $\overline{U(x)}$ are coprime in $\F_3[x]$. 
 By Lemmas \ref{NP} and \ref{binomial}, for every $i=1,2$, $N_{\ph_i}^+(F)$ has at least two sides $S_{i1}$ and $S_{i2}$ with height $1$ each. Thus $F_{S_{ij}}(y)$ is irreducible over $\F_{\ph_i}$ as it is of degree $1$ for every $i,j=1,2$. By Remark \ref{remore}, every factor $\ph_i$ provides at least two  distinct prime ideals of $\Z_K$ lying above $3$ with residue degree $f=2$ each. Thus there are at least four  distinct prime ideals of $\Z_K$ lying above $3$ with residue degree $f=2$ each. As  $x^2+1$, $x^2+x-1$, and  $x^2-x-1$ are  the unique monic irreducible polynomials of degree $2$ in $\F_3[x]$,  by Lemma \ref{index}, $3$ divides $i(K)$, and so $K$ is not monogenic.
 \item
 Similarly, if we assume that $u=1$, $m\equiv \pm 1\md{81}$, and  $v\ge 3$, then let $\ph=x^2+1$. By Lemmas \ref{NP} and \ref{binomial}, $\npp{F}$ has at least  four sides  of degree $1$ each. Thus $\ph$ provides at least four prime ideals of $\Z_K$ lying above $3$ with  residue degree $1$ each.
 \item
 Assume that $m\equiv \pm 1\md{125}$ and $t\ge 2$. Then  $\ol{F(x)}$ has two monic  factors $\ph_1$ and $\ph_2$ of degree $1$ each in $\F_5[x]$. Let $\ol{\ph}$ be one fixed factor of $\ol{F(x)}$ of degree $1$. Since $m\equiv  \pm1\md{125}$, by   Lemmas \ref{NP} and \ref{binomial}, we conclude that $\npp{F}$ has at least $3$ sides of degree $1$ each. Thus  $\ph$ provides at least $3$ prime ideals of  $\Z_K$ lying above $5$ with  residue degree $1$ each. Hence the two factors $\ph_1$ and $\ph_2$  provide   at least  $6$ prime ideals of  $\Z_K$ lying above $5$ with  residue degree $1$ each.   Therefore   $5$ divides $i(K)$, and so $K$ is not monogenic.
  \item
 Assume that $m\equiv  1\md{25}$ and $u\ge 2$. By the litle Fermat's theorem, $\ol{x^{4}-1}$ has  four distinct monic factors of degree $1$ each in $\F_5[x]$. Let $\ol{\ph}$ be one fixed factor of degree $1$ of  $\ol{F(x)}$ in $\F_5[x]$. Since $m\equiv  1\md{25}$, by   Lemmas \ref{NP} and \ref{binomial}, we conclude that $\npp{F}$ has at least $2$ sides of degree $1$ each. Thus  $\ph$ provides at least $2$ prime ideals of  $\Z_K$ lying above $5$ with  residue degree $2$ each.  It follows that  the four factors provide   at least  $8$ prime ideals of  $\Z_K$ lying above $5$ with  residue degree $1$ each.   Thus,   $5$ divides $i(K)$, and so $K$ is not monogenic. 
\item
 If $m\equiv  -1\md{25}$ and $u= 2k$,  then  $\ol{x^{12}+1}$ divides $\ol{F(x)}$ in  $\F_5[x]$ and $\ol{x^{12}+1}=\ol{(x^2+4x+2)(x^2+3x+3)(x^2+2)(x^2+x+2)(x^2+2x+3)(x^2+3)}$  in $\F_5[x]$. Let $\ol{\ph}$ be one fixed of these factors. Since $m\equiv  1\md{25}$, by   Lemmas \ref{NP} and \ref{binomial},  $\npp{F}$ has at least $2$ sides of degree $1$ each. Thus  $\ph$ provides at least $2$ prime ideals of  $\Z_K$ lying above $5$ with  residue degree $2$ each. Hence there are   at least $12$ prime ideals of  $\Z_K$ lying above $5$ with  residue degree $2$ each. As there are only $10$  monic irreducible polynomial of degree $2$ in $\F_5[x]$,  by Lemma \ref{index}, $5$ divides $i(K)$, and so $K$ is not monogenic.
 \item
 Assume that $u=v=1$, $m\equiv \pm 82\md{5^4}$, and $t\ge 3$.
Since  $\ol{x^{6}-2}=\ol{(x^2+4x+2)(x^2+x+2)(x^2+2)}$ and $\ol{x^{6}+2}=\ol{(x^2+3)(x^2+3x+3)(x^2+2x+3)}$ in $\F_5[x]$.  Let $\ol{\ph}$ be one fixed  irreducible factor of $\ol{F(x)}$ of degree $2$. Since $m\equiv \pm 82\md{5^4}$, we have  $\nu_5(m^4+1)\ge 4$. So, by   Lemmas \ref{NP} and \ref{binomial}, we conclude that $\npp{F}$ has at least $4$ sides of degree $1$ each. Thus  $\ph$ provides at least $4$ prime ideals of  $\Z_K$ lying above $5$ with  residue degree $2$ each. Since  $\ol{F(x)}$ has $3$ distinct monic irreducible factors in $\F_5[x]$ of degree $2$ each, we conclude that   there are   at least $12$ prime ideals of  $\Z_K$ lying above $5$ with  residue degree $2$ each. As there is only $10$  monic irreducible polynomial of degree $2$ in $\F_5[x]$,  by Lemma \ref{index}, $5$ divides $i(K)$, and so $K$ is not monogenic.
\end{enumerate}
\end{proof}
\begin{proof} of Theorem \ref{nspe}.\\
Since gcd$(k,30)=1$, let $(x,y)\in \z^2$ be the unique solution of the equation $k\cdot x- 2^u\cdot 3^v\cdot 5^t\cdot y=1$ and $\th=\frac{\al^x}{a^y}$. Then $\th^{2^u\cdot 3^v\cdot 5^t}=a$, and so $g(x)=x^{2^u\cdot 3^v\cdot 5^t}-a$ is the minimal polynomial of $\th$ over $\Q$; $\th\in \z_K$ is a primitive element of $K$. Since $a \neq \pm 1$ is a square free integer, we can apply Theorems \ref{pib} and \ref{npib}, and get the desired result.
\end{proof}

\end{document}